\documentclass[letterpaper, 12pt]{article}
\usepackage{amssymb,amsmath}
\usepackage{epsf}
\usepackage{times,bm,mathrsfs}
\usepackage{amsmath}
\usepackage{amssymb}
\usepackage{amsfonts}
\usepackage{mathrsfs}
\usepackage{bbm}
\usepackage{color}
\usepackage{colordvi}
\usepackage{graphicx}
\usepackage{amscd}
\usepackage{epsfig}
\newtheorem{claim}{Claim}

\definecolor{Red}{rgb}{1,0,0}
\definecolor{Blue}{rgb}{0,0,1}
\definecolor{Olive}{rgb}{0.41,0.55,0.13}
\definecolor{Green}{rgb}{0,1,0}
\definecolor{MGreen}{rgb}{0,0.8,0}
\definecolor{DGreen}{rgb}{0,0.55,0}
\definecolor{Yellow}{rgb}{1,1,0}
\definecolor{Cyan}{rgb}{0,1,1}
\definecolor{Magenta}{rgb}{1,0,1}
\definecolor{Orange}{rgb}{1,.5,0}
\definecolor{Violet}{rgb}{.5,0,.5}
\definecolor{Purple}{rgb}{.75,0,.25}
\definecolor{Brown}{rgb}{.75,.5,.25}
\definecolor{Grey}{rgb}{.5,.5,.5}
\definecolor{Black}{rgb}{0,0,0}

\newcommand{\ecal}{\mathcal{E}}
\newcommand{\fcal}{\mathcal{F}}
\newcommand{\gcal}{\mathcal{G}}

\newcommand{\real}{\mathbb{R}}

\newcommand{\eps}{\varepsilon}

\newcommand{\E}{\mathbb{E}}
\renewcommand{\P}{\mathbb{P}}

\newcommand{\glass}{ML/GLA\-SS/MT}
\newcommand{\rstar}{$R^*$/STAR/MDC}
\newcommand{\steac}{STEAC/SC}
\newcommand{\glassshort}{ML}
\newcommand{\rstarshort}{R^*}
\newcommand{\steacshort}{STEAC}
\newcommand{\frate}{\alpha}
\newcommand{\rglass}{\frate_{\mathrm{\glassshort}}}
\newcommand{\rrstar}{\frate_{\mathrm{\rstarshort}}}
\newcommand{\rsteac}{\frate_{\mathrm{\steacshort}}}

\newcommand{\phy}{S}
\newcommand{\pop}[1]{\mathrm{#1}}
\newcommand{\dt}[1]{\tau_{\pop{#1}}}
\newcommand{\dtp}[1]{\tau^{\mathrm{P}}_{\pop{#1}}}
\newcommand{\lin}[2]{I^{(#1)}_{\pop{#2}}}
\newcommand{\lout}[2]{O^{(#1)}_{\pop{#2}}}
\newcommand{\ct}[3]{T^{(#1,#3)}_{\pop{#2}}}
\newcommand{\allpops}{\mathbb{X}}
\newcommand{\ib}{t}
\newcommand{\pc}{p}
\newcommand{\fail}[1]{\textsc{Fail}_{#1}}
\renewcommand{\succ}[1]{\textsc{Success}_{#1}}
\renewcommand{\top}[1]{\mathcal{T}[#1]}
\newcommand{\ftop}[2]{Z^{(#1)}_{\mathrm{#2}}}
\newcommand{\ftoptot}[1]{\mathcal{Z}_{\mathrm{#1}}}
\newcommand{\dist}[2]{D^{(#1)}_{\mathrm{#2}}}
\newcommand{\disttot}[1]{\mathcal{D}_{\mathrm{#1}}}

\newcommand{\remove}[1]{{}}

\begin{document}

\title{An analytical comparison of
coalescent-based
multilocus methods: The three-taxon case}

\author{
Sebastien Roch\footnote{Department of Mathematics,
University of Wisconsin--Madison.}
}


\maketitle

\begin{abstract}
Incomplete lineage sorting (ILS)
is a common source of
gene tree incongruence
in multilocus analyses. 
A large number of methods have been developed
to infer species trees in the presence
of ILS.
Here we provide a mathematical analysis
of several coalescent-based
methods. Our analysis
is performed on a three-taxon species tree
and assumes that the gene trees are correctly reconstructed
along with their branch lengths. 
\end{abstract}





\section{Introduction}

Incomplete lineage sorting (ILS) is an important
confounding factor in phylogenetic analyses
based on multiple genes or loci~\cite{Maddison:97,DegnanRosenberg:09}. ILS is a 
population-level phenomenon that is
caused by the failure of two lineages 
to coalesce in a population, leading to the possibility 
that one of the lineages first coalesces 
with a lineage from a less closely related population. 
As a result, it can produce extensive gene tree incongruence
that must be accounted for appropriately in 
multilocus analyses~\cite{DegnanRosenberg:06}.

A large number of methods have been 
developed to address this source of incongruence~\cite{LiYuKu+:09}.
Several such methods rely on a statistical
model of ILS known as the multispecies coalescent. 
In this model, populations are connected by a
phylogeny. Independent coalescent processes
are performed in each population and assembled 
to produce gene trees. Several methods
have been shown to be statistically consistent under the
multispecies coalescent, that is, they are guaranteed
to return the correct species tree given enough
loci.

The performance and accuracy 
of coalescent-based multilocus methods have been 
the subject of numerous simulation studies~\cite{LeacheRannala:11,Liu+:09,YuNakhleh:10}. 
In this paper, we complement
such studies with a detailed analytical comparison
in a tractable test case, a three-taxon species tree.
We analyze 7 methods: 
maximum likelihood (ML), GLASS/Maximum Tree (MT),
$R^*$, STAR, minimizing deep coalescences (MDC),
STEAC, and shallowest coalescences (SC).
Under the assumption that gene trees are reconstructed
without estimation error, we derive the exponential
decay rate of the failure probability
as the internal branch length of the
species tree varies. The analysis,
which relies on large-deviations theory, 
reveals that 
ML and GLASS/MT 
are more accurate in this setting
than the other methods---especially in the 
regime where ILS is more common.

\section{Materials and Methods}

\subsection{Multispecies coalescent: Three-taxon case}

We first describe the statistical model under which
our analysis is performed, the {\it multispecies coalescent}.
We only discuss the three-taxon case.
For more details, see~\cite{DegnanRosenberg:09}
and references therein.

A weighted rooted tree is called ultrametric 
if each leaf is exactly at the same distance from the root.
For a three-leaf ultrametric tree $G$ with leaves
$a$, $b$, and $c$, 
we denote by $ab|c$ the topology
where $a$ and $b$ are closer to each other than to $c$, 
and similarly
for $ac|b$, $bc|a$. The topology of $G$ is denoted
by $\top{G}$.

Let $\phy$ be an ultrametric species phylogeny with three taxa. We assume that all haploid populations in $\phy$
have population size $N$. 
We denote the current
populations by $\pop{A}$,
$\pop{B}$ and $\pop{C}$ 
(which we identify with the leaves of $\phy$)
and we assume
that $\phy$ has topology
$\pop{A}\pop{B}|\pop{C}$. The ancestral
populations are $\pop{AB}$ (corresponding to the
immediate ancestor to populations $\pop{A}$ and $\pop{B}$) 
and $\pop{ABC}$ (corresponding to the ancestor
of populations $\pop{A}$, $\pop{B}$ and $\pop{C}$).
The corresponding divergence times
(backwards in time from the present)
are denoted by $\dt{AB}$ and $\dt{ABC}$ with
the assumption $\dt{AB} \leq \dt{ABC}$.
All times are given in 
units of $N$ generations.
For a population $\pop{X}$, we let $\dtp{X}$
be the divergence time of the parent population
of $\pop{X}$. Let $\allpops = \{\pop{A},\pop{B},\pop{C},\pop{AB},\pop{ABC}\}$ be the set of all populations in $S$.

We consider $L$ loci $\ell = 1,\ldots,L$ and,
for each locus, we sample
one lineage from each population at time $0$.
For locus $\ell$, we denote by $\lin{\ell}{X}$
the number of lineages entering population $\pop{X}$
and by $\lout{\ell}{X}$
the number of lineages exiting population $X$
(backwards in time), where necessarily 
$\lin{\ell}{X} \geq \lout{\ell}{X}$.
Similarly, for $k=\lout{\ell}{X}+1,\ldots,\lin{\ell}{X}$, 
the time of the coalescent event bringing the number of
lineages from $k$ to $k-1$ in population $\pop{X}$ and locus 
$\ell$ is
$\ct{\ell}{X}{k}$.
We denote by $G_1, \ldots, G_L$ the corresponding
ultrametric gene trees (including both topology and branch lengths).

Then, under the multispecies coalescent, assuming the loci are unlinked, the likelihood of
the gene trees is given by
\begin{eqnarray}
f(G_1,\ldots,G_L | \phy)
&=& \prod_{\ell=1}^L
\exp\bigg(
-\sum_{\pop{X} \in \allpops}
\bigg\{
\binom{\lout{\ell}{X}}{2}\left(\dtp{X}
- \ct{\ell}{X}{\lout{\ell}{X}+1}\right)\nonumber\\
&& \quad\quad\quad\quad\quad - \sum_{k = \lout{\ell}{X}+1}^{\lin{\ell}{X}}
\binom{k}{2} \left(\ct{\ell}{X}{k+1}
+ \ct{\ell}{X}{k}\right)\bigg\}\bigg)\label{eq:ml}
\end{eqnarray}
where we let $\ct{\ell}{X}{\lin{\ell}{X}+1} = 
\dt{X}$ for convenience~\cite{RannalaYang:03}. 

The parameter governing the extent of incomplete 
lineage sorting is the length of the internal branch
of $\phy$
$$
\ib = \dt{ABC} - \dt{AB}.
$$
The probability that the lineages from $\pop{A}$
and $\pop{B}$ fail to coalesce in branch $\pop{AB}$,
an event we denote by $\fail{\ell}$ for locus $\ell$
(and its complement by $\succ{\ell}$),
is
$$
1-\pc = e^{-\ib}.
$$
Note that, in that case, all three gene-tree topologies
are equally likely. Of course, $1 - \pc \to 1$ as $\ib \to 0$.

\subsection{Multilocus methods}

A basic goal of multilocus analyses is to reconstruct
a species phylogeny (including possibly estimates
of the divergence times) from a collection
of gene trees. Here we assume that the data consists of
$L$ gene trees $G_1, \ldots, G_L$ corresponding
to $L$ unlinked loci generated under the multispecies
coalescent. We assume further that the gene trees
are ultrametric and that their topology and 
branch lengths are estimated
without error. 

We consider several common 
multilocus methods. 
In our setting,
several of these methods are in fact equivalent
and we therefore group them below.
Note further that we only consider statistically consistent
methods, 
that is, methods that are guaranteed
to converge on the right species phylogeny as the
number of loci $L$ increases to $+\infty$ (at least,
in the test case we described above).
We briefly describe these methods. For more details, 
see e.g.~\cite{LiYuKu+:09}
and references therein.

\paragraph{\glass}
Under the multispecies coalescent,
maximum likelihood (ML) selects
the topology and divergence times
that maximizes the likelihood~\eqref{eq:ml}. 

In the GLASS method~\cite{MosselRoch:10a},
the species phylogeny is reconstructed 
from a distance matrix in which the entries are 
the minimum gene coalescence times across loci.
The equivalent Maximum Tree (MT) 
method was introduced and studied 
in~\cite{LiuPearl:07,EdLiPe:07,LiuYuPearl:10}.

A key result in~\cite{LiuYuPearl:10} is that,
in the constant-population case, the term
inside the exponential in the likelihood~\eqref{eq:ml} is 
monotonically decreasing in the divergence times.
As a result, because GLASS and MT 
select the phylogeny with the largest possible
divergence times, maximum likelihood 
is equivalent to GLASS and MT in this context.
See~\cite{LiuYuPearl:10} for details.

\paragraph{\rstar}
In the $R^*$ consensus method~\cite{Bryant:03, DeDeBr+:09}, for each three-taxon
set (here, we only have one such set), 
we include the topology that appears in highest
frequency among the loci and we reconstruct
the most resolved phylogeny that is compatible with these
three-taxon topologies.

In the STAR method~\cite{Liu+:09},
the species phylogeny is reconstructed 
from a distance matrix in which the entries are 
the average ranks of gene coalescence times across loci.
Here the root has the highest rank and the
rank decreases by one as one goes from the root to the leaves.

The minimizing deep coalescences (MDC) method~\cite{Maddison:97,ThanNakhleh:09} selects the species phylogeny
that requires the smallest number of ``extra lineages,''
that is, lineages that fail to coalesce in a
branch of the species phylogeny. 

On a three-taxon phylogeny, there is only 
three distinct rooted topologies. In each case,
the most recent divergence is assigned
rank $1$ in STAR and the other divergence
is assigned rank $2$. Hence selecting
the topology corresponding to the 
lowest average rank is equivalent 
to selecting the most common topology
among all loci---which is what $R^*$ does.
A similar argument shows that
MDC also selects the $R^*$ consensus tree
in our test case.

\paragraph{\steac}
In the STEAC method~\cite{Liu+:09},
the species phylogeny is reconstructed 
from a distance matrix in which the entries are 
the average coalescence times across loci.
The shallowest coalescences (SC) method
is similar to STEAC in that it uses average coalescence times.
The difference between the two methods is in how
they deal with multiple alleles per population. Since
we only consider the single-allele case, the two 
methods are equivalent here.

\subsection{Large-deviations approach}
\label{sec:ld}

As mentioned above, we consider estimation
methods that are statistically consistent in the sense that 
they are guaranteed
to converge on the correct species phylogeny as the
number of loci $L$ increases to $+\infty$.
To compare different methods, we derive
the rate of exponential decay of the probability
of failure. Let $S$ be a species phylogeny with
internal branch length $\ib$ and assume that
$G_1, \ldots, G_L$ are unlinked gene trees
generated under the multispecies coalescent. 
As $L \to +\infty$, 
large-deviations theory (see e.g.~\cite{Durrett:96})
gives a characterization of the {\it (exponential) decay rate}
$$
\frate_{\mathbb{M}}(\ib) 
= - \lim_{L \to +\infty} \frac{1}{L} \ln \P[\text{Method $\mathbb{M}$ fails 
given $L$ loci from $S$}].
$$
That is, roughly
$$
\P[\text{Method $\mathbb{M}$ fails 
given $L$ loci from $S$}]
\approx e^{- L \frate_{\mathbb{M}}(\ib)},
$$
for large $L$.
As the notation indicates, the key parameter that influences
the decay rate is the length of the internal branch $\ib$
of the species phylogeny. In particular, we expect that
$\frate_{\mathbb{M}}(\ib)$ is increasing in $\ib$
as a larger $\ib$ makes the reconstruction problem
easier.

To derive $\frate_{\mathbb{M}}(\ib)$,
we express the probability of failure 
as a large deviation event of the form
$$
\P[\text{Method $\mathbb{M}$ fails 
given $L$ loci from $S$}]
= \P\left[\sum_{\ell = 1}^L Y_\ell > y L\right],
$$
where $y$ is a constant and $\{Y_\ell\}_{\ell=1}^L$
are independent identically distributed 
random variables. The particular choice
of random variables depends on the method,
as we explain below. 
Let 
$$
\phi(s) = \E[e^{s Y_\ell}],
$$
be the moment-generating function of $Y_\ell$
(which does not depend on $\ell$ by assumption). 
Then the decay rate
is given by 
\begin{equation}\label{ld}
\frate_{\mathbb{M}}(\ib)
= y s_* - \ln \phi(s_*),
\end{equation}
where $s_* > 0$ is the solution (if it exists) to
$$
\frac{\phi'(s_*)}{\phi(s_*)} = y,
$$
provided there is an $s > 0$ such that
$\phi(s) < +\infty$, $y > \E[Y_\ell]$ and
$Y_\ell$ is not a point mass at $\E[Y_\ell]$.
For more details on large-deviations
theory, see e.g.~\cite{Durrett:96}.

\section{Results}
\label{sec:results}

\subsection{A domination result}
\label{sec:domination}

We first argue that, given perfectly
reconstructed unlinked gene trees under the
multispecies coalescent, \glass~always has 
a greater probability of success than
\rstar~and \steac---or, in fact, any other method. 
Indeed note that the probability of
success can be divided into two cases:
\begin{enumerate}
\item The case where 
$\succ{\ell}$ occurs for at least one locus $\ell$,
an event of probability $(1 - (1-\pc)^L)$.
In that case,
\glass~necessarily succeeds whereas
the other two methods succeed with
probability $< 1$.

\item The case where $\fail{\ell}$ occurs for
all loci $\ell$, an event of probability
$(1-\pc)^L$. In that case, all methods succeed 
with
probability $1/3$ by symmetry. For instance, 
for \glass, any pair of populations
is equally likely to lead to the smallest inter-species 
distance. A similar argument applies to the other two
methods.

\end{enumerate}
Hence, overall \glass~succeeds with greater
probability.

\subsection{Decay rates}

\begin{figure}[t]
\centering
\includegraphics[width=5in]{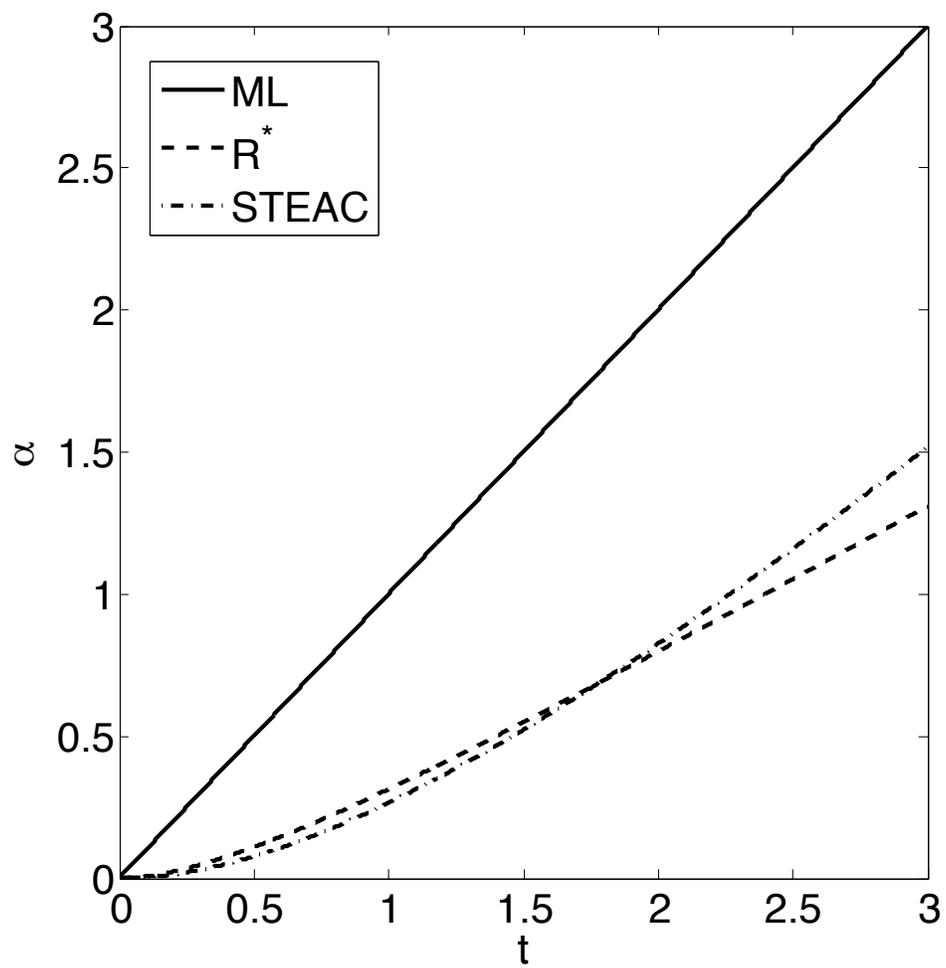}\\
\caption{
Decay rates.}\label{fig:rate1}
\end{figure}

We derive the decay rates
for the methods above. 
The results are plotted in Figure~\ref{fig:rate1}.
The asymptotic regimes are highlighted in
Figures~\ref{fig:rate01} and~\ref{fig:rate100}.
All proofs can be found
in the appendix.

\paragraph{\glass}
In this case, the decay rate can be derived
directly without using~\eqref{ld}.
Following the derivation in~\cite{MosselRoch:10a}
(see also~\cite{LiuYuPearl:10} for a similar argument),
\glass~succeeds with probability
$$
(1 - (1-\pc)^L) + \frac{1}{3}(1-\pc)^L.
$$
Then we get the following:
\begin{claim}[\glass] The decay rate
of \glass~on $S$ is
\begin{eqnarray*}
\rglass(\ib)
&=& \ib.
\end{eqnarray*}
\end{claim}

\paragraph{\rstar}
For a locus $\ell$, we let $\ftop{\ell}{AB}$
be $1$ if $\fail{\ell}$ occurs and $\top{G_\ell}
= \pop{A}\pop{B}|\pop{C}$, and $0$ otherwise.
We let
$$
\ftoptot{AB} = \sum_{\ell=1}^L \ftop{\ell}{AB}.
$$
Similarly, we define $\ftop{\ell}{AC}$,
$\ftop{\ell}{BC}$, $\ftoptot{AC}$ and $\ftoptot{BC}$.
Then \rstar~fails if
$$
\ftoptot{AB} + (L - \ftoptot{AC} - \ftoptot{BC} - \ftoptot{AB})
< \max\{\ftoptot{AC},\ftoptot{BC}\}.
$$
It can be shown that
\begin{eqnarray*}
\rrstar(\ib) &=& - \lim_{L \to +\infty} \frac{1}{L} \ln \P[
2 \ftoptot{AC} +
\ftoptot{BC}
> L].
\end{eqnarray*}
Then we get the following:
\begin{claim}[\rstar] The decay rate
of \rstar~on $S$ is
\begin{eqnarray*}
\rrstar(\ib) &=& - \ln \left(
2 \sqrt{\frac{1}{3} e^{-\ib}\left(1 - \frac{2}{3} e^{-\ib}\right)} 
+ \frac{1}{3} e^{-\ib}
\right).
\end{eqnarray*}
As $\ib \to 0$,
$$
\rrstar(\ib)
= \frac{3}{4} \ib^2 
+ O(\ib^3),
$$
and, as $\ib \to +\infty$,
$$
\rrstar(\ib) \approx \frac{\ib}{2} - \frac{1}{2}\ln \frac{4}{3}.
$$
\end{claim}

\paragraph{\steac}
For a locus $\ell$, we let $\dist{\ell}{AB}$
be the time to the most recent common ancestor
of $\pop{A}$ and $\pop{B}$ in $G_\ell$
(in units of $N$ generations).
We let
$$
\disttot{AB} = \sum_{\ell=1}^L \dist{\ell}{AB}.
$$
Similarly, we define $\dist{\ell}{AC}$,
$\dist{\ell}{BC}$, $\disttot{AC}$ and $\disttot{BC}$.
Then \steac~fails if
$$
\disttot{AB} 
> \min\{\disttot{AC},\disttot{BC}\}.
$$
It can be shown that
\begin{eqnarray*}
\rsteac(\ib)
&=&\lim_{L\to +\infty}
-\frac{1}{L} \ln \P[\disttot{AB} 
- \disttot{AC}
> 0].
\end{eqnarray*}
Then we get the following:
\begin{claim}[\steac] The decay rate
of \steac~on $S$ is
\begin{eqnarray*}
\rsteac(\ib)
&=&
- \ln \left(\frac{3e^{-s_* \ib} - s_*^2 e^{-\ib} }{3(1- s_*^2)}\right),
\end{eqnarray*}
where $0 < s_* < 1$
is the unique solution to 
the fixed-point equation
\begin{equation*}
s_*
= \frac{1}{2} [6 s_* - 3\ib(1- s_*^2)]e^{(1-s_*)\ib}.
\end{equation*}
Further, as $\ib \to 0$,
$$
\rsteac(\ib)
= \frac{3}{8} \ib^2 + O(\ib^3),
$$
and, as $\ib \to +\infty$,
$$
\rsteac(\ib) \approx \ib - \ln \ib - 0.1656.
$$
\end{claim}

\begin{figure}[t]
\centering
\includegraphics[width=5in]{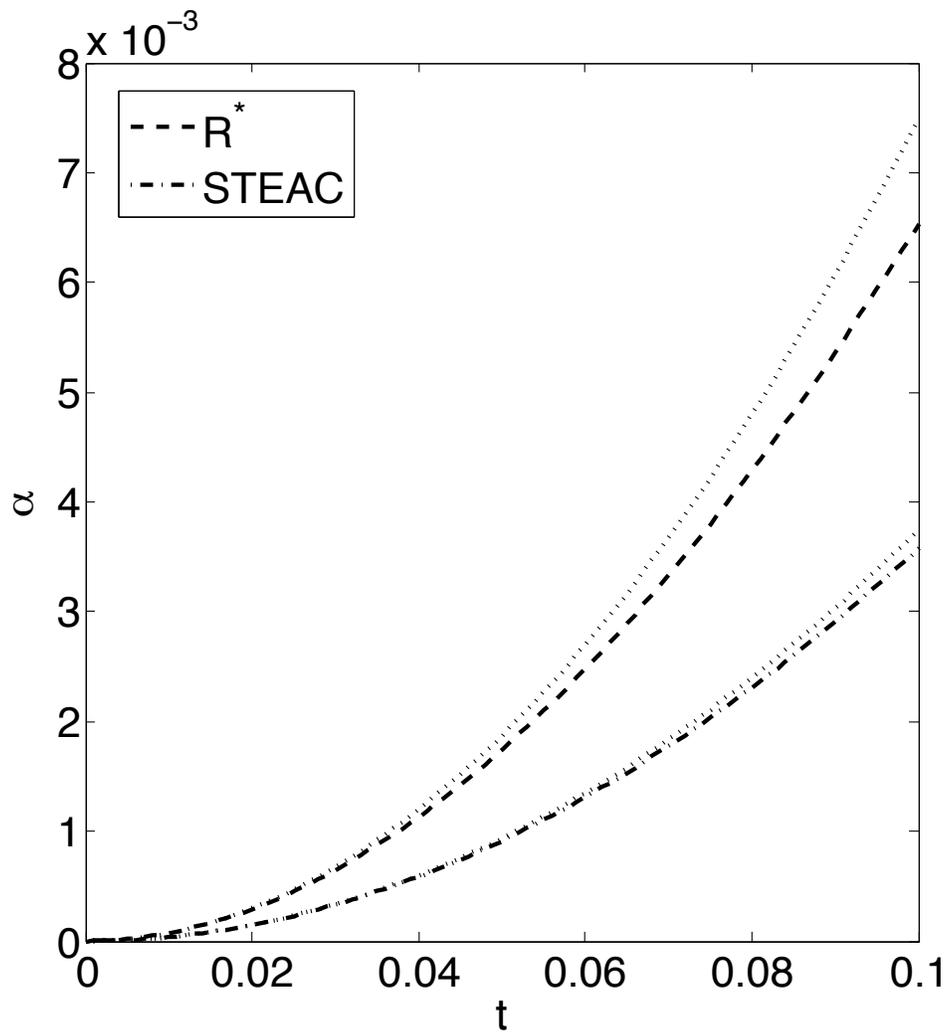}\\
\caption{
Decay rates as $\ib \to 0$. The dotted
lines indicate the respective predicted asymptotics.}\label{fig:rate01}
\end{figure}

\begin{figure}[t]
\centering
\includegraphics[width=5in]{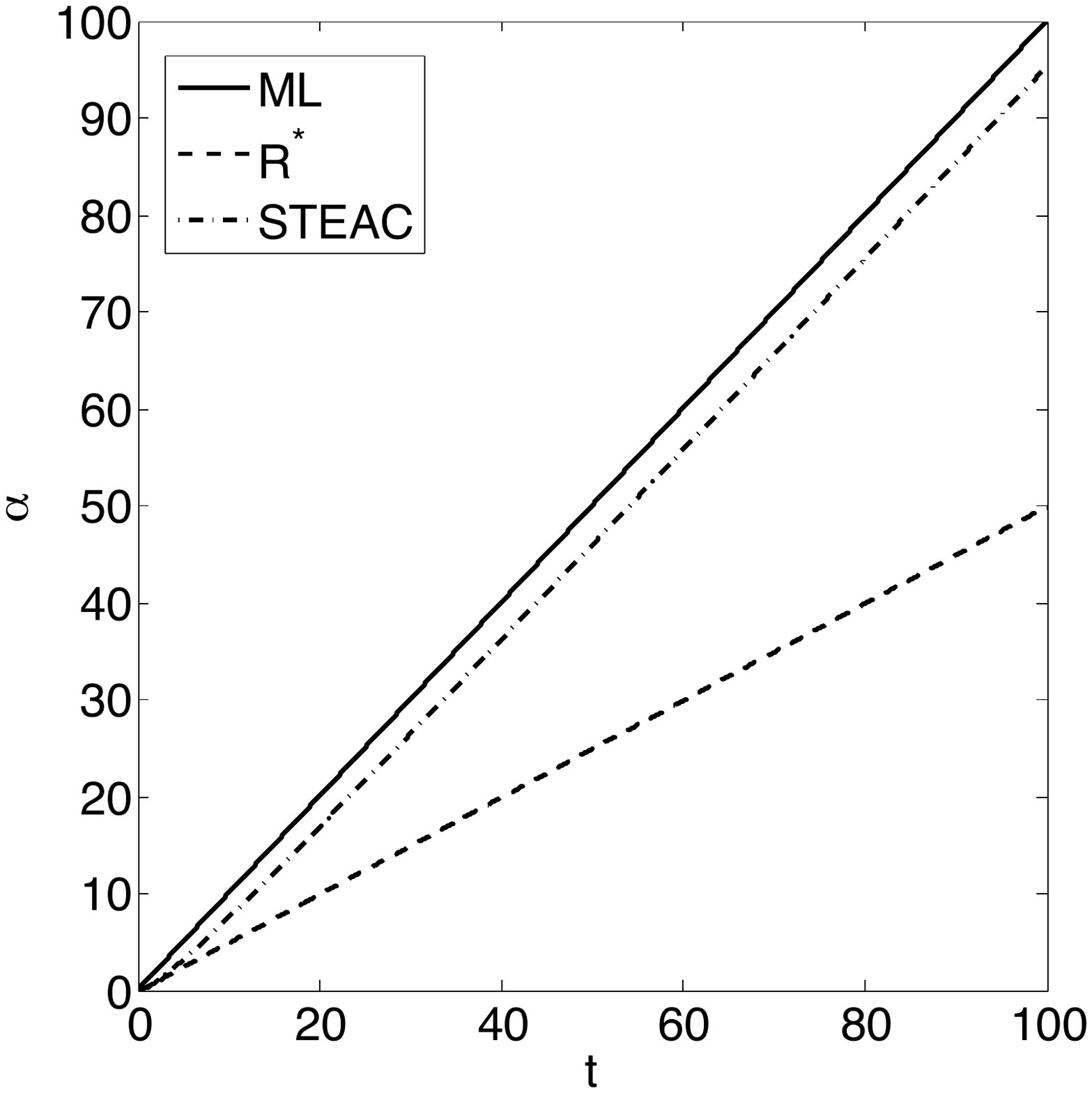}\\
\caption{
Decay rates as $\ib \to +\infty$.}
\label{fig:rate100}
\end{figure}

\section{Discussion}

As can be seen from Figures~\ref{fig:rate1},
\ref{fig:rate01} and~\ref{fig:rate100}
as well as from the asymptotics,
\glass~does indeed
give a larger decay rate for all $\ib$.
In fact, the decay rate of \glass~is
significantly higher, especially as $\ib \to 0$
that is, under high levels of incomplete lineage sorting. 
For instance,
to be concrete, if $L = 500$ loci
and $\ib = 0.1$ (in units of $N$ generations),
the probability of failure
is approximately: $1.9\times 10^{-22}$ for
\glass; $0.038$ for \rstar; $0.16$ for \steac.
Intuitively, this difference in behavior arises
from the fact that \glass~requires only
{\it one} successful locus, whereas \rstar~and \steac~rely
on an {\it average} over all loci.

Comparing \rstar~and \steac,
note that $\rrstar(\ib)$ is higher than
$\rsteac(\ib)$ for small $\ib$ but
that the situation is reversed for large $\ib$.
In fact, in the limit $\ib \to +\infty$,
$\rsteac(\ib)$ grows at roughly the same
rate as the optimal $\rglass(\ib)$.
At large $\ib$, \steac~has somewhat of an advantage
in that the expectation gap in the failure event
increases linearly with $\ib$, whereas 
it saturates under \rstar.

The analysis described here ignores
several features that 
influence the accuracy of species tree
reconstruction. Notably we have assumed
that gene trees, including their
branch lengths, are reconstructed without error.
On real sequence datasets, the uncertainty
arising from gene-tree estimation
plays an important role. For instance, although
GLASS/MT achieves the optimal decay rate
in our setting, these methods
are in fact sensitive to sequence noise
because they rely on the computation of a minimum over
loci---the very feature that leads to their
superior performance here. Extending our analysis
to incorporate gene tree estimation error is
an important open problem which should
help in the design of multilocus methods.
It is important to note that, 
under appropriate modeling of sequence
data, ML is {\it not} in general equivalent to GLASS/MT
and is likely to be more robust to
estimation error. In particular our analysis
suggest that ML may be significantly
more accurate than other methods in multilocus
studies. 

Other extensions deserve further study. 
Often many alleles are sampled
from each population. Note that the benefit of
multiple alleles is known to saturate as the
number of alleles increases~\cite{Rosenberg:02}.
This is because the probability of observing any number of alleles at the top of a branch is uniformly bounded in the number alleles existing at the bottom.  
 
Further, the molecular clock assumption,
 although it may be a reasonable first 
 approximation in the context of
 recently diverged populations, should not
be necessary for our analysis. 
One should also consider larger numbers
of taxa, varying population sizes, etc.

Simulation studies may provide further insight into 
these issues.
However an analytical approach, such as the one
we have used here, is valuable in that it allows
the study of an entire class of models in one analysis.
It can also provide useful, explicit predictions to guide
the design of reconstruction procedures.

\section{Acknowledgments}

This work was supported by NSF grant DMS-1007144
and an Alfred P. Sloan Research Fellowship.
Part of this work
was performed while the author 
was visiting the Institute for
Pure and Applied Mathematics (IPAM) at UCLA.

\clearpage

\bibliographystyle{alpha}
\bibliography{own,thesis,RECOMB12}


\clearpage

\appendix

\section{Proofs}

\subsection{\glass}

\paragraph{Decay rate}
Following the derivation in~\cite{MosselRoch:10a}
(see also~\cite{LiuYuPearl:10} for a similar argument),
\glass~succeeds with probability
$$
(1 - (1-\pc)^L) + \frac{1}{3}(1-\pc)^L,
$$
where the two terms correspond to the two cases
described in Section~\ref{sec:domination}.
Hence the decay rate of the failure probability is
\begin{eqnarray*}
\rglass(\ib) &=& \lim_{L \to +\infty} - \frac{1}{L}\ln \left\{1 - (1 - (1-\pc)^L) + \frac{1}{3}(1-\pc)^L\right\}\\
&=& \lim_{L \to +\infty} - \frac{1}{L}\ln \left\{\frac{2}{3}(1-\pc)^L\right\}\\
&=& \lim_{L \to +\infty} \left\{- \frac{1}{L}\ln\frac{2}{3}
- \ln(1-\pc)\right\}\\
&=& - \ln (1-\pc)\\
&=& \ib.
\end{eqnarray*}

\subsection{\rstar}

\paragraph{Definitions}
For a locus $\ell$, we let $\ftop{\ell}{AB}$
be $1$ if $\fail{\ell}$ occurs and $\top{G_\ell}
= \pop{A}\pop{B}|\pop{C}$, and $0$ otherwise.
We let
$$
\ftoptot{AB} = \sum_{\ell=1}^L \ftop{\ell}{AB}.
$$
Similarly, we define $\ftop{\ell}{AC}$,
$\ftop{\ell}{BC}$, $\ftoptot{AC}$ and $\ftoptot{BC}$.
Then \rstar~fails if
$$
\ftoptot{AB} + (L - \ftoptot{AC} - \ftoptot{BC} - \ftoptot{AB})
< \max\{\ftoptot{AC},\ftoptot{BC}\},
$$
an event we denote by $\ecal$.
The second term on the LHS 
comes form the fact that, given $\succ{\ell}$,
$\top{G_\ell} = \pop{A}\pop{B}|\pop{C}$.
To deal with the term on the RHS, we re-write $\ecal$ as
$$
2 \max\{\ftoptot{AC},\ftoptot{BC}\} +
\min\{\ftoptot{AC},\ftoptot{BC}\}
> L,
$$
and we use the auxiliary events 
$$
\ecal' = 
\{2 \ftoptot{AC} +
\ftoptot{BC}
> L\},
$$
and
$$
\ecal'' = 
\{2 \ftoptot{BC} +
\ftoptot{AC}
> L\},
$$
to bound $\P[\ecal]$ as follows
$$
\P[\ecal'] \leq \P[\ecal]  \leq \P[\ecal' \cup \ecal'']
\leq 2 \P[\ecal'],
$$
where we used that $\P[\ecal'] = \P[\ecal'']$ (by symmetry)
in a union bound, and the fact that, on $\ecal'$,
$$
2 \max\{\ftoptot{AC},\ftoptot{BC}\} +
\min\{\ftoptot{AC},\ftoptot{BC}\}
\geq
2 \ftoptot{AC} +
\ftoptot{BC}
> L.
$$
Hence 
$$
-\frac{1}{L} \ln \P[\ecal'] \geq 
-\frac{1}{L} \ln \P[\ecal]  \geq 
-\frac{1}{L} \ln 2 \P[\ecal']
= -\frac{1}{L} \ln \P[\ecal'] -\frac{1}{L} \ln 2,
$$
and, taking a limit as $L \to +\infty$,
$$
\rrstar(\ib) = \lim_{L \to +\infty} -\frac{1}{L} \ln \P[\ecal]
= \lim_{L \to +\infty} -\frac{1}{L} \ln \P[\ecal'],
$$
provided the limit exists.

\paragraph{Moment-generating function}
In order to compute the limit above, we use
the moment-generating function
$$
\phi(s) = \E[\exp(s[2 \ftop{\ell}{AC} + \ftop{\ell}{BC}])],
$$
(which does not depend on $\ell$)
as described in Section~\ref{sec:ld}.
Dividing up the expectation 
into the four possible cases, we have
$$
\phi(s)
= \left(\pc + \frac{1}{3}(1 - \pc)\right)
+ \frac{1}{3}(1 - \pc) (e^s + e^{2s}) < +\infty,
$$
for all $s \in \real$.
Letting 
$$
W_\pc = \frac{1}{3}(1 - \pc),
$$ 
the derivative of $\phi(s)$ is
$$
\phi'(s)
=  W_\pc (e^s + 2 e^{2s}).
$$

\paragraph{Decay rate}
By large-deviations theory,
we are looking for a solution to
$$
1 = \frac{\phi'(s)}{\phi(s)}.
$$
Letting $\omega = e^s$, we get the quadratic equation
$$
\left(\pc + W_\pc \right)
+ W_\pc (\omega + \omega^2)
= W_\pc (\omega + 2 \omega^2),
$$
or, rearranging,
$$
\left(\pc + W_\pc \right) = W_\pc \omega^2,
$$
whose solution is
$$
\omega_* = e^{s_*} = 
\sqrt{\frac{\pc + W_\pc }{W_\pc}}.
$$
Then
\begin{eqnarray*}
\rrstar(\ib)
&=& \lim_{L\to +\infty}
-\frac{1}{L} \ln \P[\ecal']\\
&=&\lim_{L\to +\infty}
-\frac{1}{L} \ln \P[2 \ftoptot{AC} +
\ftoptot{BC}
> L]\\
&=& s_*
- \ln \phi(s_*).
\end{eqnarray*}
Noting that
\begin{eqnarray*}
\phi(s_*)
&=& \pc + W_\pc 
+ W_\pc \sqrt{\frac{\pc + W_\pc }{W_\pc}}
+ W_\pc \frac{\pc + W_\pc }{W_\pc}\\
&=& 2(\pc + W_\pc)
+ \sqrt{(\pc + W_\pc)W_\pc},
\end{eqnarray*}
we get
\begin{eqnarray*}
s_*
- \ln \phi(s_*)
&=& \ln \left(
\frac{\sqrt{\pc + W_\pc}}
{2 \sqrt{W_\pc} (\pc + W_\pc)
+ W_\pc \sqrt{\pc + W_\pc}}
\right),
\end{eqnarray*}
Rearranging, we have finally
\begin{eqnarray*}
\rrstar(\ib)
&=& - \ln \left(
2 \sqrt{W_\pc(\pc + W_\pc)} 
+ W_\pc
\right)\\
&=& - \ln \left(
2 \sqrt{\frac{1}{3} e^{-\ib}\left(1 - \frac{2}{3} e^{-\ib}\right)} 
+ \frac{1}{3} e^{-\ib}
\right).
\end{eqnarray*}

\paragraph{Asymptotics}
By a Taylor expansion, we get as $\ib \to 0$ that
$$
\rrstar(\ib)
= \frac{3}{4} \ib^2 
+ O(\ib^3).
$$
On the other hand, as $\ib \to +\infty$,
\begin{eqnarray*}
\rrstar(\ib) &=& 
- \ln \left(e^{-\ib/2}
\left[2 \sqrt{\frac{1}{3}\left(1 - \frac{2}{3} e^{-\ib}\right)} 
+ \frac{1}{3} e^{-\ib/2}\right]
\right)\\
&=& \frac{\ib}{2} - \beta_\ib
\end{eqnarray*}
where
$$
\lim_{\ib \to +\infty } \beta_\ib = \frac{1}{2}\ln \frac{4}{3}. 
$$

\subsection{\steac}

\paragraph{Definitions}
For a locus $\ell$, we let $\dist{\ell}{AB}$
be the time to the most recent common ancestor
of $\pop{A}$ and $\pop{B}$ in $G_\ell$
(in units of $N$ generations).
We let
$$
\disttot{AB} = \sum_{\ell=1}^L \dist{\ell}{AB}.
$$
Similarly, we define $\dist{\ell}{AC}$,
$\dist{\ell}{BC}$, $\disttot{AC}$ and $\disttot{BC}$.
Then \steac~fails if
$$
\disttot{AB} 
> \min\{\disttot{AC},\disttot{BC}\},
$$
an event we denote by $\ecal$.
Once again, to deal with the term on the RHS, 
we re-write $\ecal$ as
$$
\disttot{AB} 
- \min\{\disttot{AC},\disttot{BC}\}
> 0,
$$
and we use the auxiliary events 
$$
\ecal' = 
\{\disttot{AB} 
- \disttot{AC}
> 0\},
$$
and
$$
\ecal'' = 
\{\disttot{AB} 
- \disttot{BC}
> 0\},
$$
to bound $\P[\ecal]$ as follows
$$
\P[\ecal'] \leq \P[\ecal]  \leq \P[\ecal' \cup \ecal'']
\leq 2 \P[\ecal'],
$$
where we used that $\P[\ecal'] = \P[\ecal'']$ (by symmetry)
in a union bound, and the fact that $\ecal'$ implies $\ecal$.
Hence 
$$
-\frac{1}{L} \ln \P[\ecal'] \geq 
-\frac{1}{L} \ln \P[\ecal]  \geq 
-\frac{1}{L} \ln 2 \P[\ecal']
= -\frac{1}{L} \ln \P[\ecal'] -\frac{1}{L} \ln 2,
$$
and, taking a limit as $L \to +\infty$,
$$
\rsteac(\ib) = \lim_{L \to +\infty} -\frac{1}{L} \ln \P[\ecal]
= \lim_{L \to +\infty} -\frac{1}{L} \ln \P[\ecal'],
$$
provided the limit exists.

\paragraph{Moment-generating function}
In order to compute the limit above, we need
the moment-generating function
$$
\phi(s) = \E[\exp(s[\dist{\ell}{AB} 
- \dist{\ell}{AC}])],
$$
(which does not depend on $\ell$).
Dividing up the expectation 
into the four possible cases, we have
\begin{eqnarray*}
\phi(s)
&=& \pc e^{-s \ib} \E\left[e^{s \tilde{E}_0}\right] \E\left[e^{-sE_0}\right]\\
&& + \frac{1}{3}(1 - \pc)  \E[e^{- s E_1}]\\
&& + \frac{1}{3}(1 - \pc)  \E[e^{s E_1}]\\
&& + \frac{1}{3}(1 - \pc)
\end{eqnarray*}
where we used:
\begin{enumerate}
\item In the case $\succ{\ell}$,  
$\dist{\ell}{AB} - \dt{AB} = \tilde{E}_0$ where
$\tilde{E}_0$ is an exponential
mean $1$ conditioned to be less than $\ib$.
Independently, using the memoryless property
of the exponential,
$\dist{\ell}{AC} - \dt{ABC} = E_0$ 
where $E_0$ is an exponential
mean $1$. Hence
$$
\dist{\ell}{AB} 
- \dist{\ell}{AC}
= \dt{AB} + \tilde{E}_0
- \dt{ABC} - E_0
= - \ib + \tilde{E}_0 - E_0.
$$

\item In the case $\fail{\ell}$ and $\top{G_\ell} = \pop{A}\pop{B}|\pop{C}$,  
$\dist{\ell}{AB} - \dt{ABC} = \tilde{E}_1$ where
$\tilde{E}_1$ the minimum of $\binom{3}{2}$
independent exponentials mean $1$, that is, 
an exponential
mean $1/ \binom{3}{2} = 1/3$.
Moreover, 
$\dist{\ell}{AC} - \dt{ABC} = \tilde{E}_1 + E_1$ 
where $E_1$ is an exponential mean $1$
independent of $\tilde{E}_1$.
Hence
$$
\dist{\ell}{AB} 
- \dist{\ell}{AC}
= \dt{ABC} + \tilde{E}_1
- \dt{ABC} - \tilde{E}_1 - E_1
= - E_1.
$$

\item In the case $\fail{\ell}$ and $\top{G_\ell} = \pop{A}\pop{C}|\pop{B}$,  
$\dist{\ell}{AC} - \dt{ABC} = \tilde{E}_1$ where
$\tilde{E}_1$ is an exponential
mean $1/ \binom{3}{2} = 1/3$.
Moreover, 
$\dist{\ell}{AB} - \dt{ABC} = \tilde{E}_1 + E_1$ 
where $\E_1$ is an exponential mean $1$
independent of $\tilde{E}_1$.
$$
\dist{\ell}{AB} 
- \dist{\ell}{AC}
= \dt{ABC} + \tilde{E}_1 + E_1
- \dt{ABC} - \tilde{E}_1
= E_1.
$$

\item In the case $\fail{\ell}$ and $\top{G_\ell} = \pop{B}\pop{C}|\pop{A}$,  
$\dist{\ell}{AC} = \dist{\ell}{AB}$. 

\end{enumerate}
Note that
$$
\E[e^{s E_0}] = \E[e^{s E_1}] = \frac{1}{1 - s},
$$
for all $|s| < 1$,
and
\begin{equation*}
\E[e^{s \tilde{E}_0}]
= \frac{1}{\pc} \int_0^\ib e^{s x} e^{-x} dx
= \frac{1 - e^{-(1-s)\ib}}{\pc(1 - s)}.
\end{equation*}
Hence
\begin{eqnarray*}
\phi(s)
&=& \pc e^{-s \ib} \frac{1 - e^{-(1-s)\ib}}{\pc(1 - s)} \frac{1}{1 + s}\\
&& + \frac{1}{3}(1 - \pc) \left(
\frac{1}{1 + s} + \frac{1}{1 - s} + 1
\right)\\
&=& \frac{e^{-s\ib} - e^{-t}}{1 - s^2}
+ \frac{1}{3} e^{-\ib} \left(\frac{3 - s^2}{1-s^2}\right)\\
&=& \frac{3e^{-s\ib} - s^2 e^{-\ib} }{3(1- s^2)}.
\end{eqnarray*}
The derivative of $\phi(s)$ is
\begin{eqnarray*}
\phi'(s)
&=& \frac{[-3\ib e^{-s\ib} - 2 s e^{-\ib} ][3(1- s^2)] - [3e^{-s\ib} - s^2 e^{-\ib} ][- 6 s]}{[3(1- s^2)]^2}\\
&=& \frac{[18 s - 9\ib(1- s^2)]e^{-s\ib} - 6 s e^{-\ib}}{[3(1- s^2)]^2}
\end{eqnarray*}

\paragraph{Decay rate}
By large-deviations theory,
we are looking for a solution to
\begin{equation}\label{eq:ldsteac}
0 = \frac{\phi'(s)}{\phi(s)} = 
\frac{[6 s - 3\ib(1- s^2)]e^{-s\ib} - 2 s e^{-\ib}}{(1- s^2)(3e^{-s\ib} - s^2 e^{-\ib})}.
\end{equation}
Note that the denominator on the RHS is positive
on $s \in (0,1)$, and that 
$$
\frac{\phi'(0)}{\phi(0)} = - \ib
$$
and
$$
\lim_{s \to 1^-} \frac{\phi'(s)}{\phi(s)} = +\infty,
$$
so that by~\cite{Durrett:96} there is a solution
$0 < s_* < 1$ to~\eqref{eq:ldsteac}. The solution $s_*$ must
satisfy
\begin{equation}\label{eq:zerosteac}
[6 s_* - 3\ib(1- s_*^2)]e^{-s_*\ib} - 2 s_* e^{-\ib} = 0.
\end{equation}
which can be re-written as 
the fixed-point equation
\begin{equation}\label{eq:fixedsteac}
s_*
= \frac{1}{2} [6 s_* - 3\ib(1- s_*^2)]e^{(1-s_*)\ib}
\equiv F_\ib(s_*),
\qquad 0 < s_* < 1.
\end{equation}
Note that $F_\ib(0) = -3\ib e^{\ib} \leq 0$ and
$F_\ib(1) = 3 > 1$.
Moreover,
\begin{eqnarray*}
F_\ib'(s)
&=& \frac{1}{2} [6 + 6\ib s] e^{(1 - s)\ib}
- \frac{\ib}{2}  [6 s - 3\ib(1- s^2)]e^{(1-s)\ib}\\
&=& \frac{1}{2}e^{(1-s)\ib} [6 + 3\ib^2(1 - s^2)] > 1,
\end{eqnarray*}
for $0 < s < 1$. 
Hence $F_\ib(s) - s$ is strictly increasing and has a
unique solution in $(0,1)$.
Eq.~\eqref{eq:fixedsteac} is easily solved numerically.

Then
\begin{eqnarray*}
\rsteac(\ib)
&=& \lim_{L\to +\infty}
-\frac{1}{L} \ln \P[\ecal']\\
&=&\lim_{L\to +\infty}
-\frac{1}{L} \ln \P[\disttot{AB} 
- \disttot{AC}
> 0]\\
&=&
- \ln \phi(s_*).
\end{eqnarray*}

\paragraph{Asymptotics}
We consider asymptotics when $\ib \to 0$. 
Define
$$
s_\eps = \frac{3}{4} \ib + \eps^{-1} \ib^2.
$$
Evaluating the LHS in~\eqref{eq:zerosteac}
(for $\eps > 0$ small but fixed) as $\ib \to 0$ gives
\begin{eqnarray*}
&&[6 s_\eps - 3\ib(1- s_\eps^2)]e^{-s_\eps\ib} - 2 s_\eps e^{-\ib}\\
&&\quad = \left[\frac{9}{2} \ib + 6 \eps^{-1} \ib^2 - 3 \ib \left(1 - \frac{9}{16} \ib^2  + O(\ib^3)\right)\right]\left[1 - \frac{3}{4}\ib^2 + O(\ib^3)\right]\\
&&\quad\quad - \left[\frac{3}{2} \ib + 2 \eps^{-1} \ib^2\right] \left[1 - \ib + \frac{\ib^2}{2} + O(\ib^3)\right]\\
&&\quad = \left[\frac{3}{2} \ib + 6 \eps^{-1} \ib^2 + O(\ib^3)\right]\left[1 - \frac{3}{4}\ib^2 + O(\ib^3)\right]\\
&&\quad\quad - \left[\frac{3}{2} \ib + 2 \eps^{-1} \ib^2\right] \left[1 - \ib + \frac{\ib^2}{2} + O(\ib^3)\right]\\
&&\quad = \left[4 \eps^{-1} + \frac{3}{2}\right]\ib^2 + O(\ib^3),
\end{eqnarray*}
so that, because
$$4 (-\eps)^{-1} + \frac{3}{2} < 0 \quad \text{and} 
\quad 4 \eps^{-1} + \frac{3}{2} > 0,$$ 
the solution of~\eqref{eq:zerosteac}
satisfies $s_{-\eps} < s_* < s_{\eps}$
for $0 < \eps < \frac{8}{3}$ and $t$ small enough.
Then
\begin{eqnarray*}
\phi(s_\eps)
&=& \frac{3e^{-s_\eps\ib} - s_\eps^2 e^{-\ib} }{3(1- s_\eps^2)}\\
&=& \frac{1}{3}
\left[3\left(1 - \frac{3}{4}\ib^2 +O(\ib^3)\right) - \left(\frac{9}{16} \ib^2 + O(\ib^3)\right) \left(1 - \ib + \frac{\ib^2}{2} + O(\ib^3)\right)\right]\\
&& \quad \times 
\left[1 + \frac{9}{16}\ib^2 + O(\ib^3)\right]\\
&=& 1 + \ib^2 \left\{- \frac{3}{4}
- \frac{3}{16} + \frac{9}{16}
\right\} + O(\ib^3).
\end{eqnarray*}
Since this holds for all $\eps > 0$ small
we get
$$
\rsteac(\ib)
=
- \ln \phi(s_*)
= \frac{3}{8} \ib^2 + O(\ib^3).
$$

For the $t \to +\infty$ asymptotics, let
$$
u = \frac{1}{t} \quad \text{and} \quad \sigma = (1 - s)t.
$$
Substituting in~\eqref{eq:fixedsteac}, we get
$$
1 - \sigma u 
= \frac{1}{2}[6(1 - \sigma u)
- 3 \sigma (2 - \sigma u)] e^{\sigma},
$$
which after rearranging
becomes
\begin{eqnarray}\label{eq:fcal}
u &=& 
\frac{3 e^\sigma - 1 - 3 \sigma e^\sigma}{3 \sigma e^\sigma -\sigma - \frac{3}{2}\sigma^2 e^\sigma}\nonumber\\
&=& 
\frac{1}{\sigma(1+ \frac{3\sigma e^\sigma}{2(3 e^\sigma - 1 - 3 \sigma e^\sigma)})}\\&\equiv& \fcal(\sigma)\nonumber.
\end{eqnarray}
We have $\fcal(0) = +\infty$. Moreover, letting
$\sigma_*$ be the only positive solution to
\begin{equation}\label{eq:sigmastar}
\gcal(\sigma_*) \equiv 3 e^{\sigma_*} - 1 - 3 \sigma_* e^{\sigma_*} = 0,
\end{equation}
we have $\fcal(\sigma_*) = 0$. Note that $\gcal'(\sigma) = - 3 \sigma e^\sigma < 0$, $\gcal(0) = 2$ and $\lim_{\sigma\to+\infty}\gcal(\sigma) = -\infty$,
so that $\sigma_*$ is well-defined.
Noticing that $\gcal$ appears in the denominator of~\eqref{eq:fcal} as well we get that $\fcal$ is strictly decreasing
between $\sigma = 0$ and $\sigma = \sigma_*$. Hence
the limit $\ib \to +\infty$ is equivalent to the limit 
$\sigma \to \sigma_*^-$. Finally, in that limit, letting $\sigma_\ib$
be the $\sigma$-value giving rise to the value $u = 1/\ib$
\begin{eqnarray*}
\rsteac(\ib) 
&=& - \ln \phi(s_*)\\
&=& - \ln \frac{3e^{-s_* \ib} - s_*^2 e^{-\ib} }{3(1- s_*^2)}\\
&=& \ib - \ln \frac{3 e^{\sigma_\ib} - (1 - \frac{\sigma_\ib}{\ib})^2}{3(1 - (1 - \frac{\sigma_\ib}{\ib})^2)}\\
&=& \ib - \ln \frac{3 e^{\sigma_\ib} - 1 + \frac{2\sigma_\ib}{\ib}
- \frac{\sigma_\ib^2}{\ib^2}}{3(\frac{2\sigma_\ib}{\ib}
- \frac{\sigma_\ib^2}{\ib^2})}\\
&=& \ib - \ln \ib - \beta_\ib,
\end{eqnarray*}
where
$$
\lim_{\ib \to +\infty} \beta_\ib 
= \ln\frac{3 e^{\sigma_*} - 1}{6 \sigma_*}
= \ln\frac{3 \sigma_* e^{\sigma_*}}{6 \sigma_*}
= \sigma_* - \ln 2,
$$
where we used~\eqref{eq:sigmastar}.

\end{document}